\documentclass[amstex,12pt, amssymb]{article}

\usepackage{mathtext}
\usepackage[cp1251]{inputenc}
\usepackage[T2A]{fontenc}
\usepackage[russian]{babel}
\usepackage[dvips]{graphicx}
\usepackage{amsmath}
\usepackage{amssymb}
\usepackage{amsxtra}
\usepackage{latexsym}
\usepackage{ifthen}

\newcounter{lemma}[section]

\newcounter{corollary}[section]

\newcounter{remark}[section]

\newcounter{theorem}[section]

\newcounter{proposition}[section]

\numberwithin{equation}{section}

\def\Xint#1{\mathchoice
   {\XXint\displaystyle\textstyle{#1}}%
   {\XXint\textstyle\scriptstyle{#1}}%
   {\XXint\scriptstyle\scriptscriptstyle{#1}}%
   {\XXint\scriptscriptstyle\scriptscriptstyle{#1}}%
   \!\int}
\def\XXint#1#2#3{{\setbox0=\hbox{$#1{#2#3}{\int}$}
     \vcenter{\hbox{$#2#3$}}\kern-.5\wd0}}
\def\dashint{\Xint-}

\def\cc{\setcounter{equation}{0}
\setcounter{figure}{0}\setcounter{table}{0}}

\begin{document}

\markboth{\centerline{Denis Kovtonyuk, Igor' Petkov and Vladimir
Ryazanov}} {\centerline {The Dirichlet problem and prime ends}}

\author{Denis Kovtonyuk, Igor' Petkov and Vladimir Ryazanov}

\title{\bf The Dirichlet problem and prime ends}

\maketitle

\large \begin{abstract} It is developed the theory of the boundary
behavior of homeomorphic solutions of the Beltrami equations
${\overline{\partial}}f=\mu\,{\partial}f$ of the Sobolev class
$W^{1,1}_{\rm loc}$ with respect to prime ends of domains. On this
basis, under certain conditions on the complex coefficient ${\mu}$,
it is proved the existence of regular solutions of its Dirichlet
problem in arbitrary simply connected domains and pseudoregular as
well as multivalent solutions in arbitrary finitely connected
domains with continuous boundary data in terms of prime
ends.\end{abstract}

\bigskip



\large \cc

\section{Introduction}
Theorems on the existence of homeomorphic solutions of the class
$W^{1,1}_{\mathrm{loc}}$ have been recently proved for many
degenerate Beltrami equations, see, e.g., the monographs \cite{GRSY}
and \cite{MRSY} and further references therein. The theory of
boundary behavior of homeomorphic solutions with generalized
derivatives and of the Dirichlet problem for the wide circle of
degenerate Beltrami equations in Jordan domains was developed in the
papers \cite{KPR}, \cite{KPRS}, \cite{KPRS_*} and \cite{RSSY}.

Recall some definitions. Let $D$ be a domain in the complex plane
${\Bbb C}$ and let $\mu:D\to{\Bbb C}$ be a measurable function with
$|\mu(z)|<1$ a.e. (almost everywhere). {\bf Beltrami equation} is an
equation of the form
\begin{equation}\label{eq1.1} f_{\bar z}=\mu(z)f_z\end{equation}
where $f_{\overline{z}}={\overline{\partial}}f=(f_x+if_y)/2$,
$f_{z}=\partial f=(f_x-if_y)/2$, $z=x+iy$, $f_x$ and $f_y$ are
partial derivatives of $f$ in $x$ and $y$, correspondingly. The
function is called its {\bf complex coefficient} and
\begin{equation}\label{eq1.1a}K_{\mu}(z)=\frac{1+|\mu(z)|}{1-|\mu (z)|}\end{equation}
its {\bf dilatation quotient}. Beltrami equation (\ref{eq1.1}) is
called {\bf degenerate} if $K_{\mu}$ is essentially unbounded, i.e.,
$K_{\mu}\notin L^{\infty}(D)$.

\medskip

Boundary values problems for the Beltrami equations are due to the
famous dissertation of Riemann who considered the partial case of
analytic functions when $\mu(z)\equiv0$, and to the works of Hilbert
(1904, 1924) and Poincare (1910) for the corresponding
Cauchy--Riemann system, see, e.g., further references in the papers
 \cite{KPRS_*} and \cite{RSSY}.

\bigskip

The classic {\bf Dirichlet problem} for Beltrami equation
(\ref{eq1.1}) in a Jordan domain $D$ is the search of continuous
function $f:D\to{\Bbb C}$, having partial derivatives of the first
order a.e., satisfying (\ref{eq1.1}) a.e. and also the boundary
condition
\begin{equation}\label{eq1002} \lim\limits_{z\to\zeta}\ {\rm
Re}\, f(z)\ =\ \varphi(\zeta)\qquad\forall\ \zeta\in\partial D
\end{equation} for a prescribed continuous function
$\varphi:\partial D\to{\Bbb R}$, see, e.g., \cite{Ve}.

\medskip

To study the similar problem in domains with more complicated
boundaries we need to apply the theory of prime ends by
Caratheodory, see, e.g., his paper \cite{Car$_2$} or Chapter 9 in
monograph \cite{CL}.

\medskip

The main difference in the case is that $\varphi$ should be already
a function of a boundary element (prime end $P$) but not of a
boundary point. Moreover, (\ref{eq1002}) should be replaced by the
condition
\begin{equation}\label{eq1002P} \lim\limits_{n\to\infty}\ {\rm
Re}\, f(z_n)\ =\ \varphi(P)
\end{equation} for all sequences of points $z_n\in D$ converging to
prime ends  $P$ of the domain $D$. Note that (\ref{eq1002P}) is
equivalent to the condition that
\begin{equation}\label{eq1002PP} \lim\limits_{z\to P}\ {\rm
Re}\, f(z)\ =\ \varphi(P)
\end{equation} along any ways in $D$ going to the prime ends $P$ of
the domain $D$.

\medskip

Later on, $\overline{D^{\prime}}_P$ denotes the completion of the
domain $D$ by its prime ends and $E_D$ denotes the space of these
prime ends, both with the topology of prime ends described in
Section 9.4 of monograph \cite{CL}. In addition, continuity of
mappings $f:\overline{D}_P\to\overline{D^{\prime}}_P$ and boundary
functions $\varphi : E_D\to\Bbb R$ should mean with respect to the
given topology of prime ends.

\medskip

Generalized homeomorphic solutions of the Beltrami equations are
mappings with finite distortion whose boundary behavior with respect
to prime ends in arbitrary finitely connected domains of $\Bbb C$
was studied in our last preprint \cite{KPR$_*$} and we refer the
reader to this text for historic comments, definitions and
preliminary remarks.

\bigskip

\cc
\section{Boundary behavior of solutions of the Beltrami equations}

On the basis of results of Section 8 in \cite{KPR$_*$}, we obtain
the corresponding results on the boundary behavior of solutions of
the Beltrami equations.

\medskip

\begin{theorem}\label{thKPR8.2F} {\it Let $D$ and $D^{\prime}$ be
bounded finitely connected domains in ${\Bbb C}$ and $f:D\to
D^{\prime}$ be a homeomorphic solution of the class $W^{1,1}_{\rm
loc}$ of (\ref{eq1.1}) with $K_{\mu}\in L^{1}(D)$. Then $f^{-1}$ is
extended to a continuous mapping of $\overline{D^{\prime}}_P$ onto
$\overline{D}_P$.}
\end{theorem}

\medskip

Furthermore, it is sufficient in Theorem \ref{thKPR8.2F} to assume
that $K_{\mu}$ is integrable only in a neighborhood of $\partial D$
or even more weak conditions which are due to Lemma 5.1 in
\cite{KPR$_*$}.

\medskip

However, any degree of integrability of $K_{\mu}$ does not guarantee
a continuous extendability of the direct mapping $f$ to the
boundary, see, e.g., an example in the proof of Proposition 6.3 in
\cite{MRSY}. Conditions for it have perfectly another nature. The
principal relevant result is the following.

\medskip

\begin{theorem}\label{t:10.1F} {\it Let $D$ and $D^{\prime}$ be
bounded finitely connected domains in ${\Bbb C}$ and $f:D\to
D^{\prime}$ be a homeomorphic solution of the class $W^{1,1}_{\rm
loc}$ of the Beltrami equation (\ref{eq1.1}) with the condition
\begin{equation}\label{ee:6.4F}\int\limits_{0}^{\varepsilon_0}
\frac{dr}{||K_{\mu}||(z_0,r)}\ =\ \infty\qquad\forall\
z_0\in\partial D\end{equation} where
$0<\varepsilon_0<d_0=\sup\limits_{z\in D}\,|z-z_0|$ and
\begin{equation}\label{eq8.7.6F} ||K_{\mu}||(z_0,r)\ =\ \int\limits_{
S({z_0},r)}K_{\mu}\ ds\ .\end{equation} Then $f$ can be extended to
a homeomorphism of $\overline{D}_P$ onto $\overline{D^{\prime}}_P$.}
\end{theorem}

\medskip

Here and later on, we set that $K_{\mu}$ is equal to zero outside of
the domain $D$.

\medskip

\begin{corollary}\label{thOSKRSS100F} {\it In particular, the
conclusion of Theorem \ref{t:10.1F} holds if
\begin{equation}\label{eqOSKRSS100dF}k_{z_0}(r)=O\left(\log{\frac1r}\right)\
\qquad\forall\ z_0\in\partial D \end{equation} as $r\to0$ where
$k_{z_0}(r)$ is the average of $K_{\mu}$ over the circle
$|z-z_0|=r$.}
\end{corollary}

\medskip

\begin{lemma}\label{lemOSKRSS1000F} {\it Let $D$ and $D^{\prime}$ be
bounded finitely connected domains in ${\Bbb C}$ and $f:D\to
D^{\prime}$ be a homeomorphic solution of the class $W^{1,1}_{\rm
loc}$ of the Beltrami equation (\ref{eq1.1}) with $K_{\mu}\in
L^{1}(D)$ and
\begin{equation}\label{eqOSKRSS1000F}
\int\limits_{\varepsilon<|z-z_0|<\varepsilon_0}K_{\mu}(z)\cdot\psi_{z_0,\varepsilon}^2(|z-z_0|)\
dm(z)\ =\ o\left(I_{z_0}^2(\varepsilon)\right)\qquad\forall\
z_0\in\partial D
\end{equation} as $\varepsilon\to0$ where
$0<\varepsilon_0<\sup\limits_{z\in D}\,|z-z_0|$ and
$\psi_{z_0,\varepsilon}(t):(0,\infty)\to[0,\infty]$,
$\varepsilon\in(0,\varepsilon_0)$, is a two-parametric family of
measurable functions such that
$$0<I_{z_0}(\varepsilon)\ :=\ \int\limits_{\varepsilon}^{\varepsilon_0}\psi_{z_0,\varepsilon}(t)\ dt\ <\ \infty
\qquad\forall\ \varepsilon\in(0,\varepsilon_0)\ .$$ Then $f$ can be
extended to a homeomorphism of $\overline{D}_P$ onto
$\overline{D^{\prime}}_P$.} \end{lemma}

\medskip

\begin{theorem}\label{thOSKRSS101F} {\it Let $D$ and $D^{\prime}$ be
bounded finitely connected domains in ${\Bbb C}$ and $f:D\to
D^{\prime}$ be a homeomorphic solution of the class $W^{1,1}_{\rm
loc}$ of the Beltrami equation (\ref{eq1.1}) with $K_{\mu}$ of
finite mean oscillation at every point $z_0\in\partial D$. Then $f$
can be extended to a homeomorphism of $\overline{D}_P$ onto
$\overline{D^{\prime}}_P$.}
\end{theorem}

\medskip

In fact, here it is sufficient for the function $K_{\mu}(z)$ to have
a dominant of finite mean oscillation in a neighborhood of every
point $z_0\in\partial D$.

\medskip

\begin{corollary}\label{corOSKRSS6.6.2F} {\it In particular, the
conclusion of Theorem \ref{thOSKRSS101F} holds if
\begin{equation}\label{eqOSKRSS6.6.3F}
\overline{\lim\limits_{\varepsilon\to0}}\ \ \
\dashint_{B(z_0,\varepsilon)}K_{\mu}(z)\ dm(z)\ <\
\infty\qquad\forall\ z_0\in\partial D\ .\end{equation}}
\end{corollary}

\begin{theorem}\label{thOSKRSS102F} {\it Let $D$ and $D^{\prime}$ be
bounded finitely connected domains in ${\Bbb C}$ and $f:D\to
D^{\prime}$ be a homeomorphic solution of the class $W^{1,1}_{\rm
loc}$ of the Beltrami equation (\ref{eq1.1}) with the condition
\begin{equation}\label{eqOSKRSS10.336aF}\int\limits_{\varepsilon<|z-z_0|<\varepsilon_0}K_{\mu}(z)\
\frac{dm(z)}{|z-z_0|^2}\ =\
o\left(\left[\log\frac{1}{\varepsilon}\right]^2\right)\qquad\forall\
z_0\in\partial D\ .\end{equation} Then $f$ can be extended to a
homeomorphism of $\overline{D}_P$ onto $\overline{D^{\prime}}_P$.}
\end{theorem}

\medskip

\begin{remark}\label{rmOSKRSS200F} Condition
(\ref{eqOSKRSS10.336aF}) can be replaced by the weaker condition
\begin{equation}\label{eqOSKRSS10.336bF}
\int\limits_{\varepsilon<|z-z_0|<\varepsilon_0}\frac{K_{\mu}(z)\,dm(z)}{\left(|z-z_0|\
\log{\frac{1}{|z-z_0|}}\right)^2}\ =\
o\left(\left[\log\log\frac{1}{\varepsilon}\right]^2\right)\qquad\forall\
z_0\in\partial D\ .\end{equation} In general, here we are able to
give a number of other conditions of logarithmic type. In
particular, condition (\ref{eqOSKRSS100dF}), thanking to Theorem
\ref{t:10.1F}, can be replaced by the weaker condition
\begin{equation}\label{eqOSKRSS10.336hF} k_{z_0}(r)\ =\ O
\left(\log\frac{1}{r}\log\,\log\frac{1}{r}\right)\ .\end{equation}
\end{remark}

Finally, we complete the series of criteria with the following
integral condition.

\begin{theorem}\label{thOSKRSS103F} {\it Let $D$ and $D^{\prime}$ be
bounded finitely connected domains in ${\Bbb C}$ and $f:D\to
D^{\prime}$ be a homeomorphic solution of the class $W^{1,1}_{\rm
loc}$ of the Beltrami equation (\ref{eq1.1}) with the condition
\begin{equation}\label{eqOSKRSS10.36bF} \int\limits_D\Phi\left(K_{\mu}(z)\right)\ dm(z)\ <\ \infty\end{equation}
for a nondecreasing convex function $\Phi:[0,\infty)\to[0,\infty)$
such that
\begin{equation}\label{eqOSKRSS10.37bF}
\int\limits_{\delta_*}^{\infty}\frac{d\tau}{\tau\Phi^{-1}(\tau)}\ =\
\infty\end{equation} at some $\delta_*>\Phi(0)$. Then $f$ can be
extended to a homeomorphism of $\overline{D}_P$ onto
$\overline{D^{\prime}}_P$.}
\end{theorem}

\medskip

\begin{corollary}\label{corOSKRSS6.6.3F} {\it In particular, the
conclusion of Theorem \ref{thOSKRSS103F} holds if at some $\alpha>0$
\begin{equation}\label{eqOSKRSS6.6.6F} \int\limits_{D}e^{\alpha
K_{\mu}(z)}\ dm(z)\ <\ \infty\ .\end{equation} }
\end{corollary}

\medskip

\begin{remark}\label{rmOSKRSS200000F} Note that condition (\ref{eqOSKRSS10.37bF})
is not only sufficient but also necessary for a continuous extension
to the boundary of all direct mappings $f$ with integral
restrictions of type (\ref{eqOSKRSS10.36bF}), see, e.g., Theorem 5.1
and Remark 5.1 in \cite{KR$_3$}. Recall also that condition
(\ref{eqOSKRSS10.37bF}) is equivalent to each of conditions
(7.14)--(7.18) in \cite{KPR$_*$}.\end{remark}

\bigskip

\section{On the Dirichlet problem in simply connected domains}

For $\varphi(P)\not\equiv{\rm const}$, $P\in E_D$, a {\bf regular
solution} of Dirichlet problem (\ref{eq1002P}) for Beltrami equation
(\ref{eq1.1}) is a continuous discrete open mapping $f:D\to{\Bbb C}$
of the Sobolev class $W_{\rm loc}^{1,1}$ with its Jacobian
\begin{equation}\label{eq1003}J_f(z)\ =\ |f_z|^2-|f_{\overline{z}}|^2\ \neq\ 0
\quad\quad{\text{a.e.}}\end{equation} satisfying (\ref{eq1.1}) a.e.
and condition (\ref{eq1002P}) for all prime ends of the domain $D$.
For $\varphi(P)\equiv c\in{\Bbb R}$, $P\in E_D$, a regular solution
of the problem is any constant function $f(z)=c+ic'$, $c'\in{\Bbb
R}$.

Recall that a mapping $f:D\to{\Bbb C}$ is called {\bf discrete} if
the pre-image  $f^{-1}(y)$ of every point $y\in{\Bbb C}$ consists of
isolated points and {\bf open} if the image of every open set
$U\subseteq D$ is open in ${\Bbb C}$. Later on, ${\Bbb D}$ denotes
the unit disk in ${\Bbb C}$.

\medskip

\begin{theorem}\label{thKPRS1001} {\it Let $D$ be a bounded simply
connected domain in $\Bbb C$ and let $\mu:{D}\to {\Bbb D}$ be a
measurable function with $K_{\mu}\in L^{1}_{\mathrm loc}$ and,
moreover,
\begin{equation}\label{eq1008}\int\limits_{0}^{\delta(z_0)}\frac{dr}{||K_{\mu}||(z_0,\,r)}\ =\ \infty
\qquad\forall\ z_{0}\in\overline{D}\end{equation} for some
$0<\delta(z_0)<d(z_0)={\sup\limits_{z\in {D}}|z-z_0|}$ and
$$||K_\mu||(z_0,\,r)\ :=\ \int\limits_{S(z_0,\,r)}K_{\mu}(z)\ ds\ .$$
Then the Beltrami equation (\ref{eq1.1}) has a regular solution $f$
of the Dirichlet problem (\ref{eq1002P}) for every continuous
function $\varphi:E_D\to{\Bbb R}$.}
\end{theorem}

\medskip

Here and later on, we set that $K_{\mu}$ is equal to zero outside of
the domain $D$.

\medskip

\begin{corollary}\label{corKPR1002} {\it In particular, the
conclusion of Theorem \ref{thKPRS1001} holds if
\begin{equation}\label{eqOSKRSS100dFD} k_{z_{0}}(\varepsilon)=O{\left(\log\frac{1}{\varepsilon}\right)}
\qquad\forall\ z_{0}\in\overline{D}\end{equation} as
$\varepsilon\to0$ where $k_{z_{0}}(\varepsilon)$ is the average of
the function $K_{\mu}$ over the circle $S(z_{0},\,\varepsilon)$.}
\end{corollary}

\medskip

{\it Proof of Theorem \ref{thKPRS1001}.} First of all note that
$E_D$ cannot consist of a single prime end. Indeed, all rays going
from a point $z_0\in D$ to $\infty$ intersect $\partial D$ because
the domain $D$ is bounded, see, e.g., Proposition 2.3 in \cite{RSal}
or Proposition 13.3 in \cite{MRSY}. Thus, $\partial D$ contains more
than one point and by the Riemann theorem, see, e.g., II.2.1 in
\cite{Goluzin}, $D$ can be mapped onto the unit disk ${\Bbb D}$ with
a conformal mapping $R$. However, then by the Caratheodory theorem
there is one-to-one correspondence between elements of $E_D$ and
points of the unit circle $\partial {\Bbb D}$, see, e.g., Theorem
9.6 in \cite{CL}.

Let $F$ be a regular homeomorphic solution of equation (\ref{eq1.1})
in the class $W_{\rm loc}^{1,1}$ which exists in view of condition
(\ref{eq1008}), see, e.g., Theorem 5.4 in paper \cite{RSY$_1$} or
Theorem 11.10 in monograph \cite{MRSY}.

Note that the domain $D^*=F(D)$ is simply connected in
$\overline{\Bbb C}$, see, e.g., Lemma 5.3 in \cite{IR} or Lemma 6.5
in \cite{MRSY}. Let us assume that $\partial D^*$ in $\overline{\Bbb
C}$ consists of the single point $\{ \infty\}$. Then $\overline{\Bbb
C}\setminus D^*$ also consists of the single point $\infty$, i.e.,
$D^*=\Bbb C$, since if there is a point $\zeta_0\in\Bbb C$ in
$\overline{\Bbb C}\setminus D^*$, then, joining it and any point
$\zeta_*\in D^*$ with a segment of a straight line, we find one more
point of $\partial D^*$ in $\Bbb C$, see, e.g., again Proposition
 2.3 in \cite{RSal} or Proposition 13.3 in \cite{MRSY}. Now, let
$\Bbb D^*$ denote the exterior of the unit disk $\Bbb D$ in $\Bbb C$
and let $\kappa(\zeta)=1/\zeta$, $\kappa(0)=\infty$,
$\kappa(\infty)=0$. Consider the mapping $F_*=\kappa\circ
F:\widetilde D\to \Bbb D_0$ where $\widetilde D=F^{-1}(\Bbb D^*)$
and $\Bbb D_0=\Bbb D\setminus \{ 0\}$ is the punctured unit disk. It
is clear that $F_*$ is also a regular homeomorphic solution of
Beltrami equation (\ref{eq1.1}) in the class $W_{\rm loc}^{1,1}$ in
the bounded two--connected domain $\widetilde D$ because the mapping
$\kappa$ is conformal. By Theorem \ref{t:10.1F} there is a
one--to--one correspondence between elements of $E_D$ and $0$.
However, it was shown above that $E_D$ cannot consists of a single
prime end. This contradiction disproves the above assumption that
$\partial D^*$ consists of a single point in $\overline{\Bbb C}$.

Thus, by the Riemann theorem $D^*$ can be mapped onto the unit disk
${\Bbb D}$ with a conformal mapping $R_*$. Note that the function
$g:=R_*\circ F$ is again a regular homeomorphic solution in the
Sobolev class $W_{\rm loc}^{1,1}$ of Beltrami equation (\ref{eq1.1})
which maps $D$ onto $\Bbb D$. By Theorem \ref{t:10.1F} the mapping
$g$ admits an extension to a ho\-meo\-mor\-phism $g_*:{\overline
D}_P\to\overline{\Bbb D}$.

We find a regular solution of the initial Dirichlet problem
(\ref{eq1002P}) in the form $f=h\circ g$ where $h$ is a holomorphic
function in $\Bbb D$ with the boundary condition
$$\lim\limits_{z\to\zeta}\,{\rm Re}\,h(z)\ =\
\varphi(g_*^{-1}(\zeta))\qquad \forall\ \zeta\in\partial{\Bbb D}\
.$$ Note that we have from the right hand side a continuous function
of the variable $\zeta$.

As known, the analytic function $h$ can be reconstructed in ${\Bbb
D}$ through its real part on the boundary up to a pure imaginary
additive constant with the Schwartz formula, see, e.g., \S\ 8,
Chapter III, Part 3 in \cite{HC},
$$h(z)\ =\ \frac{1}{2\pi i}\int\limits_{|\zeta|=1}\varphi\circ
g_*^{-1}(\zeta)\cdot\frac{\zeta+z}{\zeta-z}\cdot\frac{d\zeta}{\zeta}\
.$$ It is easy to see that the function $f=h\circ g$ is a desired
regular solution of the Dirichlet problem (\ref{eq1002P}) for
Beltrami equation (\ref{eq1.1}). $\ \Box$

\medskip

Applying Lemma 2.2 in \cite{RS}, see also Lemma 7.4 in \cite{MRSY},
we obtain the following general lemma immediately from Theorem
\ref{thKPRS1001}.

\medskip

\begin{lemma}\label{lemKPRS1000D} {\it Let $D$ be a bounded simply
connected domain in $\Bbb C$ and $\mu:{D}\to {\Bbb D}$ be a
measurable function with $K_{\mu}\in L^{1}({D})$. Suppose that, for
every $z_0\in\overline{D}$, there exist
$\varepsilon_0<d(z_0):=\sup\limits_{z\in D}|z-z_0|$ and a
one-parametric family of measurable functions
$\psi_{z_0,\,\varepsilon}:(0,\infty)\to(0,\infty)$,
$\varepsilon\in(0,\,\varepsilon_0)$ such that
\begin{equation}\label{eqKPRS1000}
0\ <\ I_{z_0}(\varepsilon)\ :=\
\int\limits_{\varepsilon}^{\varepsilon_0}
\psi_{z_0,\,\varepsilon}(t)\ dt\ < \ \infty\qquad\forall\
\varepsilon\in(0,\,\varepsilon_0)\end{equation} and as
$\varepsilon\to0$
\begin{equation}\label{eqKPRS1000a}\int\limits_{D(z_0,\,\varepsilon,\,\varepsilon_0)}
K_{\mu}(z)\cdot\psi^{2}_{z_0,\,\varepsilon}\left(|z-z_0|\right)\,
dm(z)\ =\ o(I_{z_0}^{2}(\varepsilon))\end{equation} where
$D(z_0,\,\varepsilon,\,\varepsilon_0)=\{z\in
D:\varepsilon<|z-z_0|<\varepsilon_0\}$. Then the Beltrami equation
(\ref{eq1.1}) has a regular solution $f$ of the Dirichlet problem
(\ref{eq1002P}) for every continuous function $\varphi:E_D\to{\Bbb
R}$.} \end{lemma}

\medskip

\begin{remark}\label{rmKR2.9D}
In fact, it is sufficient here to request instead of the condition
$K_{\mu}\in L^1(D)$ only a local integrability of $K_{\mu}$ in the
domain $D$ and the condition $||K_{\mu}||(z_0,r)\ne\infty$ for a.e.
$r\in (0,\varepsilon_0)$ at all $z_0\in\partial D$.
\end{remark}

\bigskip

By Lemma \ref{lemKPRS1000D} with the choice
$\psi_{z_0,\,\varepsilon}(t)\equiv 1/\left(t\log\frac{1}{t}\right)$
we obtain the following result, see also Lemma 6.1 in
\cite{KPR$_*$}.

\medskip

\begin{theorem}\label{thKPRS1000} {\it Let $D$ be a bounded simply
connected domain in $\Bbb C$ and let $\mu:{D}\to {\Bbb D}$ be a
measurable function such that
\begin{equation}\label{eq1007}{K_{\mu}(z)\leqslant Q(z)\in{\rm FMO}({\overline{D}})}\ .\end{equation}
Then the Beltrami equation (\ref{eq1.1}) has a regular solution $f$
of the Dirichlet problem (\ref{eq1002P}) for every continuous
function $\varphi:E_D\to{\Bbb R}$.} \end{theorem}

\medskip

\begin{corollary}\label{corKPR1000} {\it In particular, the
conclusion of Theorem \ref{thKPRS1000} holds if $K_{\mu}(z)\leqslant
Q(z)\in{\rm BMO}({\overline{D}})$.}
\end{corollary}

\medskip

By Corollary 6.1 in \cite{KPR$_*$} we obtain from Theorem
\ref{thKPRS1000} the following statement.

\medskip

\begin{corollary}\label{corKPR1001} {\it The conclusion of Theorem
\ref{thKPRS1000} also holds if
$$\limsup\limits_{\varepsilon\to 0}\
\dashint_{B(z_0,\,\varepsilon)}K_{\mu}(z)\ dm(z)\ <\
\infty\qquad\forall\ z_{0}\in{\overline{D}}\ .$$}
\end{corollary}

\medskip

The next statement follows from Lemma \ref{lemKPRS1000D} under the
choice $\psi(t)=1/t$, see also Remark \ref{rmKR2.9D}.

\medskip

\begin{theorem}\label{thOSKRSS102FD} {\it Let $D$ be a bounded simply
connected domain in $\Bbb C$ and let $\mu:{D}\to {\Bbb D}$ be a
measurable function such that
\begin{equation}\label{eqOSKRSS10.336aFD}\int\limits_{\varepsilon<|z-z_0|<\varepsilon_0}K_{\mu}(z)\
\frac{dm(z)}{|z-z_0|^2}\ =\
o\left(\left[\log\frac{1}{\varepsilon}\right]^2\right)\qquad\forall\
z_0\in\overline D\ .\end{equation} Then the Beltrami equation
(\ref{eq1.1}) has a regular solution $f$ of the Dirichlet problem
(\ref{eq1002P}) for every continuous function $\varphi:E_D\to{\Bbb
R}$.} \end{theorem}

\medskip

\begin{remark}\label{rmOSKRSS200FD} Similarly, choosing in Lemma
 \ref{lemKPRS1000D} $\psi(t)=1/(t\log
1/t)$ instead of $\psi(t)=1/t$, we obtain that condition
(\ref{eqOSKRSS10.336aFD}) can be replaced by the condition
\begin{equation}\label{eqOSKRSS10.336bFD}
\int\limits_{\varepsilon<|z-z_0|<\varepsilon_0}\frac{K_{\mu}(z)\,dm(z)}{\left(|z-z_0|\
\log{\frac{1}{|z-z_0|}}\right)^2}\ =\
o\left(\left[\log\log\frac{1}{\varepsilon}\right]^2\right)\qquad\forall\
z_0\in\overline D\ .\end{equation} Here we are able to give a number
of other conditions of a logarithmic type. In particular, condition
(\ref{eqOSKRSS100dFD}), thanking to Theorem \ref{thKPRS1001}, can be
replaced by the weaker condition
\begin{equation}\label{eqOSKRSS10.336hFD} k_{z_0}(r)=O
\left(\log\frac{1}{r}\log\,\log\frac{1}{r}\right).\end{equation}
\end{remark}

Finally, by Theorem \ref{thKPRS1001}, applying also Theorem 3.1 in
\cite{RSY}, we come to the following result.

\medskip

\begin{theorem}\label{thKPRS1002} {\it Let $D$ be a bounded simply
connected domain in $\Bbb C$ and let $\mu:{D}\to {\Bbb D}$ be a
measurable function such that
\begin{equation}\label{eq1009}\int\limits_{D}\Phi(K_{\mu}(z))\ dm(z)\ <\ \infty\end{equation}
where $\Phi:[0,\infty)\to[0,\infty)$ is a nondecreasing convex
function such that
\begin{equation}\label{eq1010D}\int\limits_{\delta}^{\infty}\frac{d\tau}{\tau\Phi^{-1}(\tau)}\ =\ \infty\end{equation}
for some $\delta>\Phi(0)$. Then the Beltrami equation (\ref{eq1.1})
has a regular solution $f$ of the Dirichlet problem (\ref{eq1002P})
for every continuous function $\varphi:E_D\to{\Bbb R}$.}
\end{theorem}

\medskip

\begin{remark}\label{rmkKPRS1000} Recall that condition (\ref{eq1010D})
is equivalent to each of conditions (7.14)--(7.18) in
\cite{KPR$_*$}. Moreover, condition (\ref{eq1010D}) is not only
sufficient but also ne\-ces\-sa\-ry to have a regular solution of
the Dirichlet problem (\ref{eq1002P}) for every Beltrami equation
(\ref{eq1.1}) with integral restriction (\ref{eq1009}) for every
continuous function $\varphi:E_D\to{\Bbb R}$. Indeed, by the Stoilow
theorem on representation of discrete open mappings, see, e.g.,
\cite{Sto}, every regular solution $f$ of the Dirichlet problem
(\ref{eq1002P}) for Beltrami equation (\ref{eq1.1}) with $K_{\mu}\in
L^{1}_{loc}$ can be represented in the form of composition
$f=h\circ{F}$ where $h$ is a holomorphic function and $F$ is a
regular homeomorphic solution of (\ref{eq1.1}) in the class $W_{\rm
loc}^{1,1}$. Thus, by Theorem 5.1 in \cite{RSY$_2$} on the
nonexistence of regular homeomorphic solutions of (\ref{eq1.1}) in
the class $W_{\rm loc}^{1,1}$, if (\ref{eq1010D}) fails, then there
is a measurable function $\mu:{D}\to {\Bbb D}$ satisfying integral
condition (\ref{eq1009}) for which Beltrami equation (\ref{eq1.1})
has no regular solution of the Dirichlet problem (\ref{eq1002P}) for
any nonconstant continuous function $\varphi:E_D\to{\Bbb R}$.
\end{remark}

\medskip

\begin{corollary}\label{corKPR1003} {\it In particular, the
conclusion of Theorem \ref{thKPRS1002} holds if at some $\alpha>0$}
\begin{equation}\label{eq1011}
\int\limits_{D}e^{\alpha K_{\mu}(z)}\,dm(z)<\infty\ .\end{equation}
\end{corollary}


\section{On pseudoregular solutions in multiply connected domains}

As it was first noted by Bojarski, see, e.g., \S\ 6 of Chapter 4 in
\cite{Ve}, the Dirichlet problem for the Beltrami equations,
generally speaking, has no regular solution in the class of
continuous (single--valued) in ${\Bbb C}$ functions with generalized
derivatives in the case of multiply connected domains $D$. Hence the
natural question arose: whether solutions exist in wider classes of
functions for this case ? It is turned out to be solutions for this
problem can be found in the class of functions admitting a certain
number (related with connectedness of $D$) of poles at prescribed
points. Later on, this number will take into account the
multiplicity of these poles from the Stoilow representation.


Namely, a {\bf pseudoregular solution} of such a problem is a
continuous (in $\overline{\Bbb C}$) discrete open mapping
$f:D\to\overline{\Bbb C}$ of the Sobolev class $W_{\rm loc}^{1,1}$
(outside of poles) with its Jacobian
$J_f(z)=\left|f_z\right|^2-\left|f_{\overline{z}}\right|^2\ne0$ a.e.
satisfying (\ref{eq1.1}) a.e. and the boundary condition
(\ref{eq1002P}).

\medskip

Arguing similarly to the case of simply connected domains and
applying Theorem V.6.2 in \cite{Goluzin} on conformal mappings of
finitely connected domains onto circular domains and also Theorems
4.13 and 4.14 in \cite{Ve}, we obtain the following result.

\bigskip

\begin{theorem}\label{thKPRS1001M} {\it Let $D$ be a bounded $m-$connected domain in $\Bbb
C$ with nondegenerate boundary components and $\mu:{D}\to {\Bbb D}$
be a measurable function with $K_{\mu}\in L^{1}_{\mathrm loc}$ and
\begin{equation}\label{eq1008M}\int\limits_{0}^{\delta(z_0)}\frac{dr}{||K_{\mu}||(z_0,\,r)}\ =\ \infty
\qquad\forall\ z_{0}\in\overline{D}\end{equation} for some
$0<\delta(z_0)<d(z_0)={\sup\limits_{z\in {D}}|z-z_0|}$ and
$$||K_\mu||(z_0,\,r)\ :=\ \int\limits_{S(z_0,\,r)}K_{\mu}(z)\ ds\ .$$
Then the Beltrami equation (\ref{eq1.1}) has a pseudoregular
solution $f$ of the Dirichlet problem (\ref{eq1002P}) with
$k\geqslant m-1$ poles at prescribed points in $D$ for every
continuous function $\varphi:E_D\to{\Bbb R}$.}
\end{theorem}

\medskip

Here, as before, we set $K_{\mu}$ to be extended by zero outside of
the domain $D$.

\medskip

\begin{corollary}\label{corKPR1002M} {\it In particular, the
conclusion of Theorem \ref{thKPRS1001M} holds if
\begin{equation}\label{eqOSKRSS100dFDM} k_{z_{0}}(\varepsilon)=O{\left(\log\frac{1}{\varepsilon}\right)}
\qquad\forall\ z_{0}\in\overline{D}\end{equation} as
$\varepsilon\to0$ where $k_{z_{0}}(\varepsilon)$ is the average of
the function $K_{\mu}$ over the circle $S(z_{0},\,\varepsilon)$.}
\end{corollary}

\medskip

{\it Proof of Theorem \ref{thKPRS1001M}.} Let $F$ be a regular
solution of equation (\ref{eq1.1}) in the class $W_{\rm loc}^{1,1}$
that exists by condition (\ref{eq1008M}), see, e.g., Theorem 5.4 in
the paper \cite{RSY$_1$} or Theorem 11.10 in the monograph
\cite{MRSY}.

Note that the domain $D^*=F(D)$ is $m-$connected in $\overline{\Bbb
C}$ and there is a natural one--to--one correspondence between
components $\gamma_j$ of $\gamma =\partial D$ and components
$\Gamma_j$ of $\Gamma =\partial D^*$, $\Gamma_j=C(\gamma_j, F)$ and
$\gamma_j=C(\Gamma_j, F^{-1})$, $j=1,\ldots , m$, see, e.g., Lemma
5.3 in \cite{IR} or Lemma 6.5 in \cite{MRSY}. Moreover, by Remark
1.1 in \cite{KPR$_*$} every subspace $E_j$ of $E_D$  associated with
$\gamma_j$ consists of more than one prime end, even it is
homeomorphic to the unit circle.

Next, no one of $\Gamma_j$, $j=1,\ldots , m$, is degenerated to a
single point. Indeed, let us assume that $\Gamma_{j_0}=\{ \zeta_0\}$
first for some $\zeta_0\in\Bbb C$. Let $r_0\in(0,d_0)$ where
$d_0=\inf\limits_{\zeta\in\Gamma\setminus\Gamma_{j_0}}|\zeta
-\zeta_0|$. Then the punctured disk $D_0=\{ \zeta\in\Bbb C :
0<|\zeta -\zeta_0|<r_0\}$ is in the domain $D^*$ and its boundary
does not intersect $\Gamma\setminus\Gamma_{j_0}$. Set $\widetilde
D=F^{-1}(D_0)$. Then by the construction $\widetilde D\subset D$ is
a 2--connected domain, $\overline{\widetilde
D}\cap\gamma\setminus\gamma_{j_0}=\varnothing$,
$C(\gamma_{j_0},\widetilde F)=\{ \zeta_0\}$ and
$C(\zeta_0,\widetilde F^{-1})=\gamma_{j_0}$ where $\widetilde F$ is
a restriction of the mapping $\widetilde F$ to $\widetilde D$.
However, this contradicts Theorem \ref{t:10.1F} because, as it was
noted above, $E_{j_0}$ contains more than one prime end.

Now, let assume that $\Gamma_{j_0}=\{ \infty\}$. Then the component
of $\overline{\Bbb C}\setminus D^*$ associated with $\Gamma_{j_0}$,
see Lemma 5.1 in \cite{IR} or Lemma 6.3 in \cite{MRSY}, is also
consists of the single point $\infty$ because if the interior of
this component is not empty, then choosing there an arbitrary point
$\zeta_0$ and joining it with a point $\zeta_*\in D^*$ by a segment
of a straight line we would find one more point in $\Gamma_{j_0}$,
see, e.g., Proposition 2.3 in \cite{RSal} or Proposition 13.3 in
\cite{MRSY}.

Thus, applying if it is necessary an additional stretching
(conformal mapping), with no loss of generality we may assume that
$D^*$ contains the exteriority $\Bbb D_*$ of the unit disk $\Bbb D$
in $\Bbb C$. Set $\kappa(\zeta)=1/\zeta$, $\kappa(0)=\infty$,
$\kappa(\infty)=0$. Consider the mapping $F_*=\kappa\circ
F:\widetilde D\to \Bbb D_0$ where $\widetilde D=F^{-1}(\Bbb D_*)$
and $\Bbb D_0=\Bbb D\setminus \{ 0\}$ is the punctured unit disk. It
is clear that $F_*$ is also a homeomorphic solution of Beltrami
equation (\ref{eq1.1}) of the class $W_{\rm loc}^{1,1}$ in
2--connected domain $\widetilde D$ because the mapping $\kappa$ is
conformal. Consequently, by Theorem
 \ref{t:10.1F} elements of $E_{j_0}$ should be in a one--to--one
 correspondence with $0$. However, it was already noted, $E_{j_0}$
cannot consists of a single prime end. The obtained contradiction
disproves the assumption that $\Gamma_{j_0}=\{ \infty\}$.

Thus, by Theorem V.6.2 в \cite{Goluzin}, see also Remark 1.1 in
\cite{KPR$_*$}, $D^*$ can be mapped with a conformal mapping $R_*$
onto a bounded circular domain ${\Bbb D}^*$ whose  boundary consists
of mutually disjoint circles. Note that the function $g:=R_*\circ F$
is again a regular homeomorphic solution in the Sobolev class
$W_{\rm loc}^{1,1}$ for Beltrami equation (\ref{eq1.1}) that maps
$D$ onto $\Bbb D^*$. By Theorem \ref{t:10.1F} the mapping $g$ admits
an extension to a homeomorphism $g_*:{\overline
D}_P\to\overline{\Bbb D^*}$.

Let us find a solution of the initial Dirichlet problem
(\ref{eq1002P}) in the form $f=h\circ g$ where $h$ is a meromorphic
function in $\Bbb D^*$ with the boundary condition
\begin{equation}\label{BOUNDARY}\lim\limits_{z\to\zeta}\,{\rm
Re}\,h(z)\ =\ \varphi(g_*^{-1}(\zeta))\qquad \forall\
\zeta\in\partial{\Bbb D^*} \end{equation} and $k\geqslant m-1$ poles
corresponding under the mapping $g$ to those at prescribed points in
$D$. Note that the function from the right hand side in
(\ref{BOUNDARY}) is continuous in the variable $\zeta$. Thus, such a
function $h$ exists by Theorems 4.13 and 4.14 in \cite{Ve}. It is
clear that the function $f$ associated with $h$ is by the
construction a desired pseudoregular solution of the Dirichlet
problem (\ref{eq1002P}) for Beltrami equation (\ref{eq1.1}). $\
\Box$

\medskip

Applying Lemma 2.2 in \cite{RS}, see also Lemma 7.4 in \cite{MRSY},
we obtain immediately from Theorem \ref{thKPRS1001M} the next lemma.

\medskip

\begin{lemma}\label{lemKPRS1000DM} {\it Let $D$ be a bounded $m-$connected domain in $\Bbb
C$ with nondegenerate boundary components and $\mu:{D}\to {\Bbb D}$
be a measurable function with $K_{\mu}\in L^{1}({D})$. Suppose that
for every $z_0\in\overline{D}$ there exist
$\varepsilon_0<d(z_0):=\sup\limits_{z\in D}|z-z_0|$ and
one-parameter family of measurable functions
$\psi_{z_0,\,\varepsilon}:(0,\infty)\to(0,\infty)$,
$\varepsilon\in(0,\,\varepsilon_0)$ such that
\begin{equation}\label{eqKPRS1000M}
0\ <\ I_{z_0}(\varepsilon)\ :=\
\int\limits_{\varepsilon}^{\varepsilon_0}
\psi_{z_0,\,\varepsilon}(t)\ dt\ < \ \infty\qquad\forall\
\varepsilon\in(0,\,\varepsilon_0)\end{equation} and
\begin{equation}\label{eqKPRS1000aM}\int\limits_{\varepsilon<|z-z_0|<\varepsilon_0}
K_{\mu}(z)\cdot\psi^{2}_{z_0,\,\varepsilon}\left(|z-z_0|\right)\,
dm(z)\ =\ o(I_{z_0}^{2}(\varepsilon))\ \ \ \ \ \ \  \mbox{as}\
\varepsilon\to0\ .\end{equation} Then the Beltrami equation
(\ref{eq1.1}) has a pseudoregular solution $f$ of the Dirichlet
problem (\ref{eq1002P}) with $k\geqslant m-1$ poles at prescribed
points in $D$ for every continuous function $\varphi:E_D\to{\Bbb
R}$.} \end{lemma}

\medskip

\begin{remark}\label{rmKR2.9DM}
In fact, here it is sufficient to assume instead of the condition
$K_{\mu}\in L^1(D)$ the local integrability of $K_{\mu}$ in the
domain $D$ and the condition $||K_{\mu}||(z_0,r)\ne\infty$ for a.e.
$r\in (0,\varepsilon_0)$ and all $z_0\in\partial D$.
\end{remark}

\bigskip

By Lemma \ref{lemKPRS1000DM} with the choice
$\psi_{z_0,\,\varepsilon}(t)\equiv 1/t\log\frac{1}{t}$ we obtain the
following result, see also  Lemma 6.1 in \cite{KPR$_*$}.

\medskip

\begin{theorem}\label{thKPRS1000M} {\it Let $D$ be a bounded $m-$connected domain in $\Bbb
C$ with nondegenerate boundary components and $\mu:{D}\to {\Bbb D}$
be a measurable function such that
\begin{equation}\label{eq1007M}{K_{\mu}(z)\leqslant Q(z)\in{\rm FMO}({\overline{D}})}\ .\end{equation}
Then the Beltrami equation (\ref{eq1.1}) has a pseudoregular
solution $f$ of the Dirichlet problem (\ref{eq1002P}) with
$k\geqslant m-1$ poles at prescribed points in $D$ for every
continuous function $\varphi:E_D\to{\Bbb R}$.}
\end{theorem}

\medskip

\begin{corollary}\label{corKPR1000M} {\it In particular, the
conclusion of Theorem \ref{thKPRS1000M} holds if
$K_{\mu}(z)\leqslant Q(z)\in{\rm BMO}({\overline{D}})$.}
\end{corollary}

\medskip

By Corollary 6.1 in \cite{KPR$_*$} we have by Theorem
\ref{thKPRS1000M} the next:

\medskip

\begin{corollary}\label{corKPR1001M} {\it The
conclusion of Theorem \ref{thKPRS1000M} holds if
$$\limsup\limits_{\varepsilon\to 0}\
\dashint_{B(z_0,\,\varepsilon)}K_{\mu}(z)\,dm(z)<\infty\qquad\forall\
z_{0}\in{\overline{D}}.$$} \end{corollary}

\medskip

The following statement follows from Lemma \ref{lemKPRS1000DM}
through the choice $\psi(t)=1/t$, see also Remark \ref{rmKR2.9DM}.

\medskip

\begin{theorem}\label{thOSKRSS102FDM} {\it Let $D$ be a bounded $m-$connected domain in $\Bbb
C$ with nondegenerate boundary components and $\mu:{D}\to {\Bbb D}$
be a measurable function such that
\begin{equation}\label{eqOSKRSS10.336aFDM}\int\limits_{\varepsilon<|z-z_0|<\varepsilon_0}K_{\mu}(z)\
\frac{dm(z)}{|z-z_0|^2}\ =\
o\left(\left[\log\frac{1}{\varepsilon}\right]^2\right)\qquad\forall\
z_0\in\overline D\ .\end{equation} Then the Beltrami equation
(\ref{eq1.1}) has a pseudoregular solution $f$ of the Dirichlet
problem (\ref{eq1002P}) with $k\geqslant m-1$ poles at prescribed
points in $D$ for every continuous function $\varphi:E_D\to{\Bbb
R}$.}
\end{theorem}

\medskip

\begin{remark}\label{rmOSKRSS200FDM} Similarly, choosing in Lemma
 \ref{lemKPRS1000DM} $\psi(t)=1/(t\log
1/t)$ instead of $\psi(t)=1/t$ we obtaine that the condition
(\ref{eqOSKRSS10.336aFDM}) can be replaced by the condition
\begin{equation}\label{eqOSKRSS10.336bFDM}
\int\limits_{\varepsilon<|z-z_0|<\varepsilon_0}\frac{K_{\mu}(z)\,dm(z)}{\left(|z-z_0|\
\log{\frac{1}{|z-z_0|}}\right)^2}\ =\
o\left(\left[\log\log\frac{1}{\varepsilon}\right]^2\right)\qquad\forall\
z_0\in\overline D\ .\end{equation} Here we are able to give a number
of other conditions of the logarithmic type. In particular,
condition (\ref{eqOSKRSS100dFDM}), thanking to Theorem
\ref{thKPRS1001M}, can be replaced by the weaker condition
\begin{equation}\label{eqOSKRSS10.336hFDM} k_{z_0}(r)=O
\left(\log\frac{1}{r}\log\,\log\frac{1}{r}\right).\end{equation}
\end{remark}

Finally, by Theorem \ref{thKPRS1001M}, applying also Theorem 3.1 in
the paper \cite{RSY}, we come to the following result.

\medskip

\begin{theorem}\label{thKPRS1002M} {\it Let $D$ be a bounded $m-$connected domain in $\Bbb
C$ with nondegenerate boundary components, $k\geqslant m-1$ and
$\mu:{D}\to {\Bbb D}$ be a measurable function such that
\begin{equation}\label{eq1009M}\int\limits_{D}\Phi(K_{\mu}(z))\ dm(z)\ <\ \infty\end{equation}
where $\Phi:[0,\infty)\to[0,\infty)$ is nondecreasing convex
function with the condition
\begin{equation}\label{eq1010M}\int\limits_{\delta}^{\infty}\frac{d\tau}{\tau\Phi^{-1}(\tau)}\ =\ \infty\end{equation}
for some $\delta>\Phi(0)$. Then the Beltrami equation (\ref{eq1.1})
has a pseudoregular solution $f$ of the Dirichlet problem
(\ref{eq1002P}) with $k$ poles at prescribed points in $D$ for every
continuous function $\varphi:E_D\to{\Bbb R}$.}
\end{theorem}

\medskip

Recall that condition (\ref{eq1010M}) is equivalent to every of
conditions (7.14)--(7.18) in \cite{KPR$_*$}.

\medskip

\begin{corollary}\label{corKPR1003M} {\it In particular, the
conclusion of Theorem \ref{thKPRS1002M} holds if for some
$\alpha>0$}
\begin{equation}\label{eq1011M}
\int\limits_{D}e^{\alpha K_{\mu}(z)}\,dm(z)<\infty\ .\end{equation}
\end{corollary}

\bigskip

\section{On multivalent solutions in finitely connected domains}

\medskip

In multiply connected domains ${D}\subset{\Bbb{C}}$, in addition to
pseudoregular solutions, Dirichlet problem (\ref{eq1002P}) for
Beltrami equations (\ref{eq1.1}) admits multivalent solutions in the
spirit of the theory of multivalent analytic functions.

\medskip

We say that a discrete open mapping
$f:B(z_0,\varepsilon_0)\to{\Bbb{C}}$, where
$B(z_0,\varepsilon_0)\subset{D}$, is {\bf a local regular solution}
of equation (\ref{eq1.1}) if $f\in{W_{\rm loc}^{1,1}}$,
$J_{f}\neq{0}$ and $f$ satisfies (\ref{eq1.1}) a.e. Two local
regular solutions $f_0:B(z_0,\,\varepsilon_0)\to{\Bbb{C}}$ and
$f_*:B(z_*,\,\varepsilon_*)\to{\Bbb{C}}$ of equation (\ref{eq1.1})
is called an {\bf extension of each to other} if there is a chain of
such solutions $f_i:B(z_i,\varepsilon_i)\to{\Bbb{C}}$,
$i=\overline{1,m}$, that $f_1=f_0$, $f_m=f_*$ and
$f_{i}(z)\equiv{f_{i+1}(z)}$ for
${z\in{E_i}}:=B(z_i,\,\varepsilon_{i})\cap{B(z_{i+1},\,\varepsilon_{i+1})}\neq\varnothing$,
$i=\overline{1,m-1}$. A collection of local regular solutions
$f_{j}:B(z_{j},\varepsilon_{j})\to{\Bbb{C}}$, $j\in J$, is said to
be a {\bf multivalent solution of equation} (\ref{eq1.1}) in $D$, if
the disks $B(z_j,\,\varepsilon_j)$ cover the whole domain $D$ and
$f_j$ are mutually extended each to other through this collection
and the collection is maximal by inclusion. A multivalent solution
of (\ref{eq1.1}) is called {\bf multivalent solution of Dirichlet
problem} (\ref{eq1002P}) if $u(z)=Re{f(z)}=Re{f_j(z)}$, $z\in
B(z_j,\,\varepsilon_j)$, $j\in J$, is a single--valued function in
$D$ and satisfies condition (\ref{eq1002P}).

\medskip

The proof of the existence of multivalent solutions of Dirichlet
problem (\ref{eq1002P}) for Beltrami equation (\ref{eq1.1}) in
finitely connected domains is reduced to the Dirichlet problem for
harmonic functions in circular domains, see, e.g., \S\ 3 of Chapter
VI in \cite{Goluzin}.

\medskip

\begin{theorem}\label{thKPRS3001} {\it Let $D$ be a bounded $m-$connected domain in $\Bbb
C$ with nondegenerate boundary components and $\mu:{D}\to {\Bbb D}$
be a measurable function which satisfies hypotheses of Theorems
\ref{thKPRS1001M}--\ref{thKPRS1002M} or Corollaries
\ref{corKPR1002M}--\ref{corKPR1003M}. Then Beltrami equation
(\ref{eq1.1}) has a multivalent solution of Dirichlet problem
(\ref{eq1002P}) for every continuous function $\varphi:E_{D}\to{\Bbb
R}$.} \end{theorem}

{\it Proof.} It is sufficient to prove the statement of the theorem
under the hypotheses of Theorem \ref{thKPRS1001M} because the
hypotheses of the rest theorems and corollaries imply the hypotheses
of Theorem \ref{thKPRS1001M} as it was shown above.

Next, similarly to the first part of Theorem \ref{thKPRS1001M}, we
first prove that there is a regular homeomorphic solution $g$ of
Beltrami equation (\ref{eq1.1}) mapping the domain $D$ onto a
circular domain $\Bbb D^*$ whose boundary consists of mutually
disjoint circles. By Theorem \ref{t:10.1F} the mapping $g$ admits an
extension to a homeomorphism $g_*:{\overline D}_P\to\overline{\Bbb
D^*}$.

As known, in the circular domain $\Bbb D^*$, there is a solution of
the Dirichlet problem
\begin{equation}\label{BOUNDARYM}\lim\limits_{z\to\zeta}\ u(z)\ =\ \varphi(g_*^{-1}(\zeta))\qquad \forall\
\zeta\in\partial{\Bbb D^*} \end{equation} for harmonic functions
$u$, see, e.g., \S\ 3 of Chapter VI in \cite{Goluzin}. Let
$B_0=B(z_0,r_0)$ is a disk in the domain $D$. Then ${\cal B}_0 =
g(B_0)$ is a simply connected subdomain of the circular domain $\Bbb
D^*$ where there is a conjugate function $v$ determined up to an
additive constant such that $h=u+iv$ is a single--valued analytic
function. The function $h$ can be extended to, generally speaking
multivalent, analytic function $H$ along any path in $\Bbb D^*$
because $u$ is given in the whole domain  $\Bbb D^*$.

Thus, $f=H\circ g$ is a desired multivalent solution of the
Dirichlet problem (\ref{eq1002P}) for Beltrami equation
(\ref{eq1.1}).  $\ \Box$

\bigskip

\begin{remark}\label{rmkKPRS3001} Note that it can be proved an analog of the
known theorem on monodromy for analytic functions stating that any
multivalent solution of Beltrami equation (\ref{eq1.1}) in a simply
connected domain $D$ is its regular single--valued solution.

Note also, here the hypothesis that the boundary components of the
domain is not degenerate to single points is essential as it is
shown by the simplest case $\mu(z)\equiv0$ of analytic functions in
the punctured unit disk because the isolated singularities of
harmonic functions are removable and by the maximum principle
harmonic functions in the unit disk are uniquely determined by its
continuous boundary values.
\end{remark}

\medskip
\noindent
{\bf Denis Kovtonyuk, Igor' Petkov and\\ Vladimir Ryazanov,}\\
Institute of Applied Mathematics and Mechanics,\\
National Academy of Sciences of Ukraine,\\
74 Roze Luxemburg Str., Donetsk, 83114, Ukraine,\\
denis$\underline{\ \ }$\,kovtonyuk@bk.ru, igorpetkov@list.ru,\\
vl$\underline{\ \ }$\,ryazanov@mail.ru

\end{document}